\documentclass[12pt]{amsart}

\usepackage{amsmath}
\usepackage{amssymb}

\newtheorem{thm}{Theorem}%[section]
\newtheorem{lem}[thm]{Lemma}%[section]
\theoremstyle{definition}
\newtheorem{defn}{Definition}%[section]

\theoremstyle{remark}
\theoremstyle{plain}

\numberwithin{equation}{section}

\def\CC{{\mathbb C}}

\def\NN{{\mathbb N}}
\def\QQ{{\mathbb Q}}
\def\RR{{\mathbb R}}

\def\ZZ{{\mathbb Z}}

\def\H{{\mathfrak H}}

\def\e{\operatorname{e}}
\def\i{\operatorname{i}}

\def\C{\operatorname{C}}

\def\L{\operatorname{L}}

\def\SL{\operatorname{SL}}

\def\lcm{\operatorname{lcm}}

\def\sgn{\operatorname{sgn}}

\def\scrB{{\mathcal B}}

\def\scrH{{\mathcal H}}

\def\scrM{{\mathcal M}}

\def\Re{\operatorname{Re}}
\def\Im{\operatorname{Im}}

\begin{document}

\title{Holomorphic almost modular forms}
\author{Jens Marklof}
\address{School of Mathematics, University of Bristol,
Bristol BS8 1TW, U.K.} 
\address{{\tt j.marklof@bristol.ac.uk}}
\thanks{January 2003/September 2003. To appear in the Bulletin of the LMS}

\begin{abstract}
Holomorphic almost modular forms 
are holomorphic functions of the complex upper half plane
which can be approximated arbitrarily well (in a suitable sense)
by modular forms of congruence subgroups of large index in $\SL(2,\ZZ)$.
It is proved that such functions have a 
rotation-invariant limit distribution when the argument approaches
the real axis.
An example for a holomorphic almost modular form is
the logarithm of $\prod_{n=1}^\infty (1-\exp(2\pi\i n^2 z))$.
The paper is motivated by the author's studies [J.\,Marklof,
Int. Math. Res. Not. {\bf 39} (2003) 2131-2151] on 
the connection between almost modular functions and
the distribution of the sequence $n^2x$ modulo one.
\end{abstract}

\maketitle

\section{Introduction}

{\em Almost modular functions} have recently been introduced in
connection with the distribution of the sequence $n^2 x$ modulo one
($n=1,2,3,\ldots$). It is shown in \cite{almost} that,
for every piecewise smooth periodic function $\psi:\RR\to\CC$ of period one, 
and for $x$ uniformly distributed in $[0,1)$, the error term
\begin{equation}\label{R}
R^x_\psi(M)= \frac{1}{\sqrt M}
\bigg( \sum_{n=1}^M \psi(n^2 x)-M \int_0^1 \psi(t)\, dt \bigg)
\end{equation}
has a limit distribution as $M\to\infty$, which can be identified with
the value distribution of a certain almost modular function.
This observation resembles results by Heath-Brown \cite{Heath92}
and Bleher \cite{Bleher92}, \cite{Bleher99}, who prove that error
terms in lattice point problems for convex planar domains
have limit distributions associated with {\em almost periodic functions}
(in the sense of Besicovitch).

In the present work we draw attention to the
holomorphic species of almost modular functions. 
A {\em holomorphic almost modular form} ({\sc hamf}) is defined 
as a holomorphic function of the complex upper half plane $\H$ to $\CC$
that can be approximated arbitrarily well 
(in a sense to be made precise in Sect.~\ref{secalmost})
by modular forms of congruence subgroups $\Gamma_1(N)$ in $\SL(2,\ZZ)$, 
as $N\to\infty$.
An example for a {\sc hamf} is
the logarithm of
\begin{equation}\label{ex}
\prod_{n=1}^\infty \big(1-e(n^2 z)\big),
\end{equation}
where $e(z):=\exp(2\pi\i z)$, cf.~Sect.~\ref{secalmost}.
As a consequence of the general limit theorem for almost modular functions
\cite[Th.~8.2]{almost},
we will see in Sect.~\ref{secLimit} that,
for $\Re z$ uniformly distributed in $[0,1)$,
\begin{equation}\label{ex2}
(\Im z)^{1/4} \log \prod_{n=1}^\infty \big(1-e(n^2 z)\big)
\end{equation}
has a rotation-invariant limit distribution in $\CC$, 
as $\Im z\to 0$. An explicit formula for
the variance of the limit distribution is given in Sect.~\ref{secVar}.
The paper concludes with a short appendix  (Section \ref{app})
containing background material from \cite{almost} 
on the definition of general almost modular functions and their limit theorems.

\section{Holomorphic modular forms}

For any integer $x$ and any prime $p$ the standard quadratic residue symbol 
$\big(\frac{x}{p}\big)$ is $1$ if $x$ is a square modulo $p$,
and $-1$ otherwise.
The {\em generalized quadratic residue symbol}
$\big(\frac{a}{b}\big)$ is,
for any integer $a$ and any odd integer $b$,
characterized by the properties \cite[pp.~160-161]{Lion80}, 
\begin{itemize}
\item[(i)] $\big(\frac{a}{b}\big)=0$ if $(a,b)\neq 1$.
\item[(ii)] $\big(\frac{a}{-1}\big)=\sgn a$.
\item[(iii)] If $b>0$, $b=\prod_i b_i$, $b_j$ primes, not necessarily distinct,
then $\big(\frac{a}{b}\big)=\prod_i \big(\frac{a}{b_i}\big)$.
\item[(iv)] $\big(\frac{a}{-b}\big)=\big(\frac{a}{-1}\big)
\big(\frac{a}{b}\big)$.
\item[(v)] $\big(\frac{0}{\pm 1}\big)= 1$.
\end{itemize}
It follows from these properties that the symbol is bimultiplicative:
\begin{equation}
\left(\frac{a_1a_2}{b}\right)
=\left(\frac{a_1}{b}\right)\left(\frac{a_2}{b}\right),
\qquad
\left(\frac{a}{b_1 b_2}\right)
=\left(\frac{a}{b_1}\right)\left(\frac{a}{b_2}\right) .
\end{equation}
Furthermore, if $b>0$, then $\big(\frac{\cdot}{b}\big)$ defines a character
modulo $b$; 
if $a\neq 0$, then $\big(\frac{a}{\cdot}\big)$ defines a character
modulo $4a$.

The action of 
$\gamma=\left(\begin{smallmatrix} a & b \\ c & d \end{smallmatrix}\right)
\in\SL(2,\RR)$ on the complex upper half plane
$\H=\{ z\in\CC : \Im z>0 \}$ is defined by the fractional linear
transformation $z\mapsto \frac{az+b}{cz+d}$.
We are here interested in the congruence subgroups
\begin{equation}
\Gamma_1(N)=\left\{
\begin{pmatrix}
a & b \\ c & d
\end{pmatrix}\in\SL(2,\ZZ) \;:\; a\equiv d\equiv 1,\,
c \equiv 0 \bmod N \right\} .
\end{equation}
Fix $\Gamma=\Gamma_1(N)$, with $4|N$.

A {\em holomorphic modular form of weight $\kappa$ for $\Gamma$} 
(with $\kappa\in\frac12\ZZ$) 
is a holomorphic function $\H\rightarrow\CC$ that satisfies 
the functional relation
\begin{equation}\label{holo}
f(\gamma z)=  \left(\frac{c}{d}\right)^{2\kappa} (c z + d)^\kappa f(z)
\end{equation}
for all
$\gamma=\left(\begin{smallmatrix} a & b \\ c & d \end{smallmatrix}\right)
\in\Gamma$,
and that is holomorphic with respect to each cusp 
\cite[Def.~1.3.3.]{Sarnak90}.
This means that $f$ has a Fourier expansion of the form
\begin{equation}
\sum_{m=0}^\infty \widehat f_m^{(i)} \, e(m z_i)  
\end{equation}
for each cuspidal coordinate $z_i$ (cf.~the appendix, Section \ref{app});
hence $f$ is bounded in each cusp.
In (\ref{holo}) the square root 
$z^{1/2}$ is chosen such that $-\pi/2<\arg z^{1/2}\leq\pi/2$,
and $z^{m/2}:=(z^{1/2})^m$, for $m\in\ZZ$.

Famous examples of holomorphic modular forms are the theta series
\begin{equation}
\theta(z)=\sum_{n\in\ZZ} e(n^2 z),
\end{equation}
which is of weight $\kappa=\frac12$, and Jacobi's $\Delta$-function
\begin{equation}
\Delta(z)= e(z) \prod_{n=1}^\infty \big(1-e(nz)\big)^{24},
\end{equation}
where $\kappa=12$, see \cite{Sarnak90}.

\begin{lem}\label{modform}
Let $a_1,a_2,\ldots,a_K\in\CC$. Then the function
\begin{equation}
\xi^{({K})}(z)= \sum_{k=1}^K a_k \, \theta(kz)
\end{equation}
is a modular form of weight $\frac12$ for $\Gamma_1(N)$
with $N=4\lcm(2,3,\ldots,K)$.
\end{lem}

\begin{proof}
We have for all 
$\gamma=\left(\begin{smallmatrix} a & b \\ c & d \end{smallmatrix}\right)
\in\Gamma_1(4k)$
\begin{equation}
\theta\left(k \frac{az+b}{cz+d}\right)
=
\theta\left(\frac{a(k z)+k b}{(c/k)(k z)+d}\right)
=
\left(\frac{(c/k)}{d}\right) (c z + d)^{1/2} \theta(kz),
\end{equation}
because 
\begin{equation}
\begin{pmatrix}
a & kb \\
c/k & d
\end{pmatrix}
\in\Gamma_1(4).
\end{equation}
Since the generalized quadratic residue symbol is multiplicative,
\begin{equation}
\left(\frac{c}{d}\right) = \left(\frac{(c/k)}{d}\right) 
\left(\frac{k}{d}\right).
\end{equation}
Furthermore $\left(\frac{k}{\,\cdot\,}\right)$ is a character mod $4k$,
and hence, for $d\equiv 1\bmod 4k$, we have
\begin{equation}
\left(\frac{k}{d}\right) =\left(\frac{k}{1}\right) = 1.
\end{equation}
This shows that $\theta^{(k)}(z):=\theta(kz)$ is a modular form
for $\Gamma_1(4k)$. The lemma now follows from the observation that
\begin{equation}
\Gamma_1(N) \subset \bigcap_{k=1}^K \Gamma_1(4k) .
\end{equation}
\end{proof}

\section{Holomorphic almost modular forms}\label{secalmost}

\begin{defn}
We call a holomorphic periodic function
$\H\to\CC$ 
\begin{equation}
\xi(z)=\sum_{m=0}^\infty \widehat \xi_m \, e(mz)
\end{equation}
a {\em holomorphic almost modular form {\sc (hamf)} of weight $\frac12$} 
if for every $\epsilon>0$ there is an $N=N(\epsilon)$ and a 
modular form 
\begin{equation}
f_\epsilon(z)= \sum_{m=0}^\infty \widehat f_{\epsilon,m} \, e(mz)
\end{equation} 
of weight $\frac12$ for $\Gamma_1(N)$, 
such that 
\begin{equation}\label{a3}
\limsup_{M\rightarrow\infty} \frac{1}{\sqrt{M}} \sum_{m=0}^M 
\big| \widehat\xi_m - \widehat f_{\epsilon,m} \big|^2 < \epsilon^2.
\end{equation}
\end{defn}

To construct examples of such functions, let $h:\H \to \CC$
be a periodic holomorphic function of the form
\begin{equation} \label{h1}
h(z)= \sum_{k=1}^\infty \widehat h_k \, e(kz),
\end{equation}
with constants $C>0$ and $\beta>\frac14$ such that for all $k\in\NN$
\begin{equation} \label{h2}
| \widehat h_k | \leq \frac{C}{k^\beta} .
\end{equation}

\begin{thm} \label{main}
For $h$ as in {\rm (\ref{h1})} and {\rm (\ref{h2})},
the function
\begin{equation}
\xi(z)=\sum_{n=1}^\infty h(n^2 z)
\end{equation}
is a {\sc hamf} of weight $\frac12$.
\end{thm}

\begin{proof}
We choose as approximants the modular forms (cf.~Lemma \ref{modform})
\begin{equation}
\xi^{(K)}(z):=\frac12 \sum_{k=1}^K \widehat h_k \, \theta(kz) .
\end{equation}
The Fourier coefficients of $\xi(z)$ are
\begin{equation}
\widehat\xi_0=0, \qquad
\widehat\xi_m
=\sum_{k=1}^\infty \sum_{\substack{n=1 \\ m=n^2 k}}^\infty \widehat h_k  ,
\end{equation}
and those of $\xi^{(K)}(z)$
\begin{equation}
\widehat\xi^{(K)}_0=\frac{1}{2} \sum_{k=1}^K \widehat h_k , \qquad
\widehat\xi^{(K)}_m
=\sum_{k=1}^K \sum_{\substack{n=1 \\ m=n^2 k}}^\infty \widehat h_k .
\end{equation}
Therefore
\begin{multline}\label{deri}
\sum_{m=1}^M \left| \widehat\xi_m-\widehat\xi^{(K)}_m\right|^2
= \sum_{k_1,k_2=K+1}^\infty 
\sum_{\substack{n_1,n_2=1 \\ 1\leq n_1^2 k_1=n_2^2 k_2\leq M}}^\infty  
\widehat h_{k_1} \overline{\widehat h_{k_2}}  \\
= \sum_{(p,q,r,s)\in S_1} 
\widehat h_{r p^2} \overline{\widehat h_{r q^2}}  
\leq C^2 \sum_{(p,q,r,s)\in S_1} 
\frac{1}{r^{2\beta} p^{2\beta} q^{2\beta}} 
\end{multline}
where the sums are restricted to the set
\begin{equation}
S_1=\{ p,q,r,s\in\NN, \quad \gcd(p,q)=1, \quad  rp^2,rq^2 > K, \quad
1\leq r p^2 s^2 q^2 \leq M \} .
\end{equation}
Thus
\begin{equation}\label{deri2}
\sum_{m=1}^M \left| \widehat\xi_m-\widehat\xi^{(K)}_m\right|^2
\leq C^2 \sqrt{M} \sum_{(p,q,r)\in S_2}
\frac{1}{r^{\frac12+2\beta} p^{1+2\beta} q^{1+2\beta}} + O(1) 
\sum_{(p,q,r)\in S_3}
\frac{1}{r^{2\beta} p^{2\beta} q^{2\beta}},
\end{equation}
where
\begin{equation}
S_2=\{ p,q,r\in\NN, \quad \gcd(p,q)=1, \quad  rp^2,rq^2 > K \} ,
\end{equation}
\begin{equation}
S_3=\{ p,q,r\in\NN, \quad \gcd(p,q)=1, \quad
1\leq r p^2 q^2 \leq M \} .
\end{equation}
The last sum in (\ref{deri2}) is bounded by (assume without loss of generality
that $\frac14<\beta<\frac12$)
\begin{equation}
\sum_{(p,q,r)\in S_3}
\frac{1}{r^{2\beta} p^{2\beta} q^{2\beta}}
=
O(M^{1-2\beta})
\sum_{\substack{p,q=1  \\
1\leq p^2 q^2 \leq M}}^\infty
\frac{1}{p^{2-2\beta} q^{2-2\beta}} 
= O(M^{1-2\beta}) .
\end{equation}
Hence
\begin{equation}\label{rem}
\sum_{m=1}^M \left| \widehat\xi_m-\widehat\xi^{(K)}_m\right|^2
\leq C^2 \sqrt{M} \sum_{(p,q,r)\in S_2}
\frac{1}{r^{\frac12+2\beta} p^{1+2\beta} q^{1+2\beta}}
+ O(M^{1-2\beta}).
\end{equation}
For $\beta>\frac14$, the sum in (\ref{rem}) converges and tends to zero for
$K$ large.
The remainder in (\ref{rem}) is of sub-leading order, 
i.e., $O(M^{1-2\beta})=o(\sqrt M)$.
So, given any $\epsilon>0$, there is a large $K$ such that
\begin{equation}
\limsup_{M\to\infty}
\frac{1}{\sqrt{M}} 
\sum_{m=0}^M \left| \widehat\xi_m-\widehat\xi^{(K)}_m\right|^2
< \epsilon^2 .
\end{equation}
\end{proof}

If we choose $h(z)=\log(1-e(z))$, we have $\widehat h_k=-1/k$
and thus Theorem \ref{main} implies that
\begin{equation}
\xi(z)=\log \prod_{n=1}^\infty \big(1-e(n^2 z)\big)
\end{equation}
is a {\sc hamf}.

\section{The limit theorem}\label{secLimit}

\begin{thm}\label{limit}
Let $\xi(z)$ be a {\sc hamf} of weight $\frac12$.
Then, for $x:=\Re z$ uniformly distributed in $[0,1)$
with respect to Lebesgue measure, 
$(\Im z)^{1/4}\xi(z)$ has a limit distribution as $y:=\Im z\rightarrow 0$.
That is, there exists a probability measure $\nu_\xi$ on $\CC$ such that,
for any bounded continuous function $g:\CC\rightarrow\CC$, we have
\begin{equation}
\lim_{y\rightarrow 0} \int_0^1 g\big(y^{1/4} \xi(x+\i y) \big) dx
= \int_\CC g(w)\, \nu_\xi(dw).
\end{equation}
Furthermore $\nu_\xi$ is invariant under rotations about the origin.
\end{thm}

\begin{proof}
The following two lemmas show that Theorem \ref{limit}
is a special case of the limit theorem for almost modular functions,
Theorem \ref{limitthm}. Rotational invariance 
of the limit distribution is proved at the end of this section.
\end{proof}

The manifold $\scrM_N=\Delta_1(N)\backslash\widetilde\SL(2,\RR)$
and the function spaces $B_\sigma(\scrM_N)$,
$\scrB^2$, which appear below, are defined in the appendix  
(Section \ref{app}).

\begin{lem}
If $f(z)$ is a holomorphic modular form of weight $\kappa$ for $\Gamma_1(N)$,
then the function
\begin{equation}
F(z,\phi)= (\Im z)^{\kappa/2} f(z) \e^{-\i \kappa \phi}
\end{equation}
is a modular function of class $B_{\kappa/2}(\scrM_N)$.
\end{lem}

\begin{proof}
For $[\gamma,\beta_\gamma]\in\Delta_1(N)$ we have
\begin{equation}
F\big([\gamma,\beta_\gamma](z,\phi)\big)
=\big[\Im(\gamma z)\big]^{\kappa/2} 
f(\gamma z) \e^{-\i \kappa (\phi+\beta_\gamma)}
=F(z,\phi)
\end{equation}
because of (\ref{holo}) and 
\begin{equation}
\big[\Im(\gamma z)\big]^{\kappa/2}=\frac{(\Im z)^{\kappa/2}}{|cz+d|^\kappa}, \qquad
\e^{-\i \kappa \beta_\gamma} = j_\gamma(z)^{-2\kappa}
=\left\{ \left(\frac{c}{d}\right) 
\left(\frac{cz+d}{|cz+d|}\right)^{1/2} \right\}^{-2\kappa}.
\end{equation}
Hence $F$ is a smooth function on $\scrM_N$.
Because (by definition) $f$ is bounded in each cusp, we have 
\begin{equation}
F(z,\phi)= O(y_i^{\kappa/2}),
\end{equation}
cf.~condition (\ref{growth}) in the appendix (Section \ref{app}).
\end{proof}

\begin{lem}
If $\xi(z)$ is a {\sc hamf} of weight $\frac12$,
then the function
\begin{equation}
\Xi(z)= (\Im z)^{1/4} \xi(z)
\end{equation}
is an almost modular function of class $\scrB^2$.
\end{lem}

\begin{proof}
If $f_\epsilon$ are the approximants of $\xi$, we choose
as approximants for $\Xi$ the functions
\begin{equation}
F_\epsilon(z,\phi)=(\Im z)^{1/4} f_\epsilon(z) \e^{-\i \phi/2}.
\end{equation}
Then, in view of (\ref{a3}),
\begin{equation}
\limsup_{y\rightarrow 0} \int_0^1 
\big| \Xi(x+\i y) - F_\epsilon(x+\i y,0) \big|^2 dx
=\limsup_{y\rightarrow 0}  y^{1/2} 
\sum_{m=0}^\infty \big| \widehat\xi_m - \widehat f_{\epsilon,m} \big|^2
\e^{-4\pi m y}
=O(\epsilon^2) ,
\end{equation}
and hence $\Xi\in\scrB^2$, see Definition \ref{alm1}.
\end{proof}

\begin{proof}[Proof of rotational invariance]
The limit distribution for every approximant 
$f_\epsilon$ is given by (cf.~Theorem \ref{horoequi})
\begin{equation}
\int_\CC g(w)\, \nu_{f_\epsilon}(dw)
= \int_{\scrM_N}
g\big(y^{1/4} f_\epsilon(x + \i y) \e^{-\i\phi/2}\big) 
\frac{dx\,dy\,d\phi}{y^2}.
\end{equation}
Substituting $\phi+2\omega$ for $\phi$
shows that $f_\epsilon(z)$ and $f_\epsilon(z)\e^{-\i \omega}$
have the same limit distribution for all $\omega\in[0,2\pi)$.
Hence $\xi(z)$ and $\xi(z)\e^{-\i \omega}$ share the same
limit distribution.
\end{proof}

\section{The variance}\label{secVar}

Let us return to the example introduced in Theorem \ref{main}
and derive an explicit formula for the variance of the
limit distribution.
By slightly modifying the steps in (\ref{deri}), (\ref{deri2}) one finds that
for 
\begin{equation}
\alpha(t):=\sum_{0\leq m <t } \big| \widehat\xi_m \big|^2
\end{equation}
and $t\to\infty$ we have
\begin{equation}
\alpha(t)
\sim A\, t^{1/2} , \qquad \text{with }
A:= \sum_{r=1}^\infty 
\sum_{\substack{p,q=1 \\ \gcd(p,q)=1 }}^\infty
\frac{\widehat h_{rp^2} \overline{\widehat h_{r q^2}}}{pq\sqrt{r}} .
\end{equation}
Therefore, using Parseval's equality and a 
standard Abelian Theorem for the Laplace transform 
\cite[Chap.~V]{Widder41},
\begin{equation}
\int_0^1 | \xi(x+\i y) |^2 dx = 
\sum_{m=0}^\infty \big| \widehat\xi_m \big|^2 \e^{-4\pi m y}=
\int_0^\infty \e^{-4\pi y t} d\alpha(t)
\sim A\, \Gamma(\tfrac32)\, (4\pi y)^{-1/2},
\end{equation}
as $y\to 0$. Since Euler's function evaluates to $\sqrt\pi/2$, we have
\begin{equation}
\lim_{y\to 0}
y^{1/2} \int_0^1 | \xi(x+\i y) |^2 dx = 
\frac14 \sum_{r=1}^\infty 
\sum_{\substack{p,q=1 \\ \gcd(p,q)=1 }}^\infty
\frac{\widehat h_{rp^2} \overline{\widehat h_{r q^2}}}{pq\sqrt{r}} . 
\end{equation}

\section{Appendix: Almost modular functions for $\widetilde\SL(2,\RR)$
\label{app}}

This appendix provides some background material on the
definition of almost modular functions and their limit theorems;
the reader is referred to \cite{almost} for more detailed information.

Let us denote by $\C(\H)$ the space of continuous
functions $\H\rightarrow\CC$, and put $\epsilon_g(z)=(cz+d)/|cz+d|$. 
The universal covering group  of $\SL(2,\RR)$ is defined as the set
\begin{equation}
\widetilde\SL(2,\RR) = \{ [g,\beta_g]\; :\; g\in\SL(2,\RR), \;
\beta_g\in\C(\H) \text{ such that } 
\e^{\i\beta_g(z)}=\epsilon_g(z) \} , 
\end{equation}
with multiplication law 
\begin{equation} \label{law}
[g,\beta^1_g][h, \beta^2_h] = [gh, \beta^3_{gh}], 
\quad \beta^3_{gh}(z)=\beta^1_g(hz)+\beta^2_h(z).
\end{equation}
We may identify $\widetilde\SL(2,\RR)$ with $\H\times\RR$ via
$[g,\beta_g]\mapsto (z,\phi)=(g\i,\beta_g(\i))$. 
The action of $\widetilde\SL(2,\RR)$ on $\H\times\RR$ 
is then $[g,\beta_g] (z,\phi) = (gz, \phi+\beta_g(z))$.
The Haar measure of $\widetilde\SL(2,\RR)$ reads, in this parametrization,
\begin{equation}\label{haar}
d\mu(g)= \frac{dx\, dy\, d\phi}{y^2}.
\end{equation}

The group $\Delta_1(N)$ is the following discrete subgroup of  
$\widetilde\SL(2,\RR)$,
\begin{equation}
\Delta_1(N)=\left\{ [\gamma,\beta_\gamma]\; :\; \gamma\in\Gamma_1(N),
\;
\beta_\gamma\in\C(\H) \text{ such that }  
\e^{\i\beta_\gamma(z)/2}= j_\gamma(z) \right\} ,
\end{equation}
where 
\begin{equation}\label{totter}
j_\gamma(z)=\left(\frac{c}{d}\right) \left(\frac{cz+d}{|cz+d|}\right)^{1/2},
\qquad
\gamma=
\begin{pmatrix}
a & b \\ c & d
\end{pmatrix}\in\Gamma_1(4).
\end{equation}
Here $z^{1/2}$ denotes the principal branch of the square-root of $z$, i.e.,
the one for which $-\pi/2<\arg z^{1/2}\leq\pi/2$.

The homogeneous space $\scrM_N=\Delta_1(N)\backslash\widetilde\SL(2,\RR)$
has finite volume with respect to Haar measure (\ref{haar}).
$\scrM_N$ has a finite number of cusps, which are represented
by the set $\eta_1,\ldots,\eta_\kappa\in\QQ \cup \infty$
on the boundary of $\H$. 
Let $\gamma_i$ be a fractional linear 
transformation $\H\to\H$ which maps the cusp at $\eta_i$ to 
the standard cusp at $\infty$ of width one. Thus 
$(z_i,\phi_i)=\widetilde\gamma_i (z,\phi)$ yields a new set of coordinates,
where the $i$th cusp appears as a cusp at $\infty$, which is
invariant under $(z_i,\phi_i)\mapsto(z_i+1,\phi_i)$.
The variable $y_i=\Im(\gamma_i z)$ measures the height into the $i$th cusp.
For any $\sigma\geq 0$, we denote by $B_\sigma(\scrM_N)$ the class of 
functions $F\in\C(\scrM_N)$ such that, for all $i=1,\ldots,\kappa$,
\begin{equation}\label{growth}
F(z,\phi)=O(y_i^\sigma), \qquad y_i\rightarrow\infty,
\end{equation}
where the implied constant is independent of $(z,\phi)$.
In view of the form of the invariant measure (\ref{haar}) we note 
that $B_\sigma(\scrM_N)\subset \L^p(\scrM_N,\mu)$ if $\sigma<1/p$.

The following theorem \cite[Th.~6.1]{almost} states that
the closed horocycles
\begin{equation}
\Delta_1(N)\{ (x+\i y,0): x\in[0,1)\}
\end{equation} 
are asymptotically equidistributed in $\scrM_N$, as $y\to 0$.

\begin{thm} \label{horoequi}
Let $0\leq\sigma<1$. Then, for every $F\in B_\sigma(\scrM_N)$,
we have 
\begin{equation}
\lim_{y\rightarrow 0} \int_0^1 F(x+\i y, 0) \, dx 
= \frac{1}{\mu(\scrM_N)} \int_{\scrM_N} F\, d\mu .
\end{equation}
\end{thm}

Let us now turn to the definition of {\em almost modular functions
of class $\scrB^p$} or $\scrH$, respectively, as given in
\cite[Sect.~7]{almost}. In the following we will consider functions 
$\Xi:\H\rightarrow\CC$, which are {\em periodic}, i.e., for 
which $\Xi(z+1)=\Xi(z)$.

\begin{defn}\label{alm1}
For any $p\geq 1$, let $\scrB^p$ 
be the class of periodic functions
$\Xi:\H\rightarrow\CC$ with the property that for every
$\epsilon>0$ there are an integer $N=N(\epsilon)>0$ and a function 
$F_\epsilon\in  B_\sigma(\scrM_N)$ with $0\leq\sigma<1/p$ so that 
\begin{equation}\label{a1}
\limsup_{y\rightarrow 0} \int_0^1 
\big| \Xi(x+\i y) - F_\epsilon(x+\i y,0) \big|^p dx
< \epsilon^p .
\end{equation}
\end{defn}

\begin{defn}\label{alm2}
Let $\scrH$ be the class of periodic functions
$\Xi:\H\rightarrow\CC$ 
with the property that for every
$\epsilon>0$ there are an integer $N=N(\epsilon)>0$ and a 
bounded continuous function $F_\epsilon\in\C(\scrM_N)$ 
such that
\begin{equation}\label{a2}
\limsup_{y\rightarrow 0} \int_0^1 
\min\big\{ 1 , \big| \Xi(x+\i y) - F_\epsilon(x+\i y,0) \big|\big\} dx
< \epsilon .
\end{equation}
\end{defn}

If $1\leq q\leq p$ we have the inclusion
$\scrB^p \subset \scrB^q  \subset \scrH$, see \cite[Prop.~7.3]{almost}.
The central observation of \cite{almost} is the following limit theorem
for almost modular functions \cite[Th.~8.2]{almost}.

\begin{thm} \label{limitthm}
Let $\Xi\in\scrH$.
Then, for $x$ uniformly distributed in $[0,1)$ with respect to
Lebesgue measure, 
$\Xi(x+\i y)$ has a limit distribution as $y\rightarrow 0$.
That is, there exists a probability measure $\nu_\Xi$ on $\CC$ such that,
for every bounded continuous function $g:\CC\rightarrow\CC$, 
\begin{equation}
\lim_{y\rightarrow 0} \int_0^1 g\big(\Xi(x+\i y)\big) dx
= \int_\CC g(w)\, \nu_\Xi(dw). 
\end{equation}
\end{thm}

The proof of this theorem follows closely the argument 
for almost periodic functions \cite{Bleher92}. The main difference
is that the equidistribution theorem for irrational Kronecker flows on 
multidimensional tori is here replaced by Theorem \ref{horoequi},
cf.~\cite[Sect.~8]{almost}.

\subsection*{Acknowledgements}
This research has been supported by an EPSRC Advanced
Research Fellowship, the Nuffield Foundation (Grant NAL/00351/G) and
the EC Research Training Network (Mathematical Aspects of Quantum Chaos) 
HPRN-CT-2000-00103.

\end{document}